\documentclass[12pt]{amsart}

\usepackage{amsmath}
\usepackage{amsfonts}
\usepackage{amssymb}
\usepackage{algorithm,algorithmic}
\usepackage{multirow}

\usepackage[usenames]{color}

\usepackage{hyperref}
\usepackage{cleveref}

\def\g{\mathfrak{g}}
\def\G{\mathfrak{G}}
\def\h{\mathfrak{h}}

\def\F{\mathbb{F}}
\def\R{\mathbb{R}}
\def\N{\mathbb{N}}
\def\Z{\mathbb{Z}}

\def\ad{\operatorname{ad}}

\def\Der{\operatorname{Der}}

\def\dim{\operatorname{dim}}
\def\End{\operatorname{End}}

\def\Ker{\operatorname{Ker}}

\def\Id{\operatorname{Id}}

\def\Rad{\operatorname{Rad}}

\newtheorem{Theorem}{Theorem}[section]
\newtheorem{Lemma}[Theorem]{Lemma}
\newtheorem{Prop}[Theorem]{Proposition}

\newtheorem{Cor}[Theorem]{Corollary}

\newtheorem{Remark}[Theorem]{Remark}

\author[Rodr\'iguez-Vallarte et-al]
{M.C. Rodr\'iguez-Vallarte${}^{(1)}$,
G. Salgado${}^{(1)}$,
O.A. S\'anchez-Valenzuela${}^{(2)}$}

\title[On extension of F-K and S Lie algebras]
{On Extensions of Frobenius-K\"ahler and Sasakian Lie Algebras}

   \address{(1) Facultad de Ciencias, UASLP; Av. Parque Chapultepec 1570, 
   Priv. del Pedregal, CP 78210, San Luis Potos\'i, SLP, M\'exico.}
   
   \email{mcvallarte@fc.uaslp.mx}
   \email{gsalgado@fciencias.uaslp.mx, gil.salgado@gmail.com}
   
   \address{(2) Centro de Investigaci\'on en Matem\'aticas, A.C., Unidad M\'erida, Yucat\'an, 
   CP 97302, M\'exico}
   \email{adolfo@cimat.mx}


   
  \date{\today}
   
\begin{document}

   
\begin{abstract}
Extensions of Lie algebras equipped with Sasakian or Frobenius-K\"ahler
geometric structures are studied.
Conditions are given so that a double extension of a Sasakian Lie algebra
be Sasakian again. 
Conditions are also given for obtaining either a Sasakian or a Frobernius-K\"ahler 
Lie algebra upon respectively extending a Frobernius-K\"ahler or a Sasakian Lie algebra
by adjoining a derivation of the source algebra.
Low-dimensional examples are included.
\end{abstract}
   
   \maketitle
   
   \tableofcontents
   

\section{Introduction}

\medskip
There are two well known mechanisms to extend a given Lie algebra $\g$.
One is through central extensions using $2$-cocycles 
$\theta\in(\wedge^2\g^*)\otimes V$ with values in some $\g$-module $V$.
The other is through semidirect sums
by adjoining to $\g$ a Lie subalgebra of its Lie algebra of derivations 
$\Der(\g)$. Both of them have been thoroughly studied and they coincide
when the $2$-cocycle is trivial and the adjoined Lie subalgebra
of $\Der(\g)$ is abelian and acts trivially on $\g$.
In both cases one ends up with the direct sum of two Lie algebras.

\medskip
When the two procedures are applied successively, the result is
known as a {\bf double extension} (of $\g$) and it was originally
introduced by V. Ka\c{c} in \cite{Kac}. The simplest example of 
this process starts with a Lie algebra $\g$. Then perform 
a central extension on it via a $2$-cocycle
$\theta$ with values in the underlying field viewed as the 
trivial $1$-dimensional $\g$-module. Then adjoin to the
resulting central extension, a derivation of the central extension itself.

\medskip
The relevance of this process lies in the fact that,
when applied to a Lie algebra equipped with some kind of geometric structure
(e.g., an invariant metric, a symplectic structure, a hermitian structure, 
a contact structure, etc.), the resulting Lie algebra inherits
the same kind of geometric structure.
See \cite{Kac}, \cite{D-M-1}, \cite{M-R-2}, \cite{RV-S}, \cite{RV-S-2}.

\medskip
As a matter of fact, the double extension process has been extended
from the realm of Lie algebras to almost any kind of algebraic structure;
i.e., Leibniz algebras, Poisson algebras, K\"ahler Lie algebras; Hom-Lie algebras,
and their $\Z_2$-graded versions; say, Lie superalgebras; Leibniz superalgebras; 
Hom-Lie superalgebras; etc.
In other words, the double extension technique has
proved to be useful and effective. 
See \cite{A-B}, \cite{B-B}, \cite{B-H}, \cite{B-M}, \cite{M-R}.

\medskip
For example, it can be proved that if $\g$ comes 
equipped with a symplectic form $\omega:\g\times\g\to\R$
---a fact referred to by saying that $(\g,\omega)$
{\bf is a symplectic Lie algebra}---
the central extension of $\g$ by the $2$-cocycle defined by $\omega$
is a contact Lie algebra.
Similarly, it has been proved in \cite{A-F-V} that the
central extension of a K\"ahler Lie algebra through the
$2$-cocycle of its symplectic structure, is a Sasakian Lie algebra.
It is also known that any contact Lie algebra with non-trivial center
is a central extension of a symplectic Lie algebra (see \cite{Diatta}). 
On the other hand, it is proved 
in \cite{A-F-V} that any Sasakian Lie algebra
with non-trivial center is a central extension of a K\"ahler Lie algebra.
Unfortunately, it has not been possible so far to prove
similar results starting with either contact Lie algebras or Sasakian
Lie algebras (see  Proposition \ref{Prop:KahlerNoCentralExtensionSasakian}).
 It has been proved in \cite{B-R-S}, however, that it is possible
to go from contact Lie algebras to Frobenius Lie algebras (i.e.,
symplectic Lie algebras whose symplectic form is exact) and
viceversa, while in \cite{RV-S} conditions are given so that 
a double extension of a contact Lie algebra be again
a contact Lie algebra. It is then natural to try and see
whether these kind of results can also be obtained for 
either the class of K\"ahler Lie algebras or for the class
of Sasakian Lie algebras.

\medskip
Since double extensions of K\"ahler Lie algebras have already
been studied (see \cite{D-M} and \cite{D-M-1}), the purpose of
this work is to focus on the conditions needed on a double
extension of a Sasakian Lie algebra, to be Sasakian again. 
We shall also prove that certain subclass of K\"ahler Lie algebras,
once extended by a derivation, always yield Sasakian Lie algebras
and viceversa.

\medskip
This work is organized as follows: We first recall in \S2 the basic
definitions of K\"ahler and Sasakian Lie algebras, together with
the main results of the double extension process. It is
proved in \S3 a Theorem that gives the conditions under which a double extension
of a Sasakian Lie algebra is Sasakian again.
Section \S4 deals with extensions of Frobenius-K\"ahler Lie algebras
and of Sasakian Lie algebras by a derivation, and the relationship
between them. Finally, \S5 gives a few non-trivial low-dimensional
examples of the results we have obtained in \S3 and \S4.


\section{Frobenius-K\"ahler and Sasakian Lie Algebras}


\subsection{Basic definitions} 

Let $\g$ be a Lie algebra with Lie bracket $[\,\cdot\,,\,\cdot\,]:\g\times\g\to\g$.
Any $2$-cocycle $\theta\in\wedge^2\g^*$ defines a 
{\bf central extension of} $\g$ in the following way:
If $z\notin\g$, consider the vector space
$\g_{\theta}:= \g \oplus \langle z \rangle$, and define
$[\,\cdot\,,\,\cdot\,]_{\theta}:\g_{\theta}\times\g_{\theta}\to\g_{\theta}$
by means of 
$[x,y]_{\theta} = [x,y]_{\g} + \theta (x,y) z$, for any $x$ and $y$ in $\g$,
and $[x,z]_{\theta}=0$, for any $x$ in $\g$.
Since $\theta\in\wedge^2\g^*$ is a $2$-cocycle, $[\,\cdot\,,\,\cdot\,]_{\theta}$
is in fact a Lie bracket on $\g_{\theta}$. The pair $(\g_{\theta},[\,\cdot\,,\,\cdot\,]_{\theta})$
is called {\bf the central extension of $\g$ by the $2$-cocycle $\theta$}.
Moreover, it is proved in \cite{C-E} that  $\g_{\theta_1}\simeq\g_{\theta_2}$
if and only if the cohomology classes $[\theta_1]$ and $[\theta_2]$ 
are the same in $H^2(\g,\F)$, where $\F$ is the base field.

\medskip
One may also consider the case in which $\g$ is extended by
adjoining a single derivation $D\in \Der(\g)$. Such an extension
has the vector space decomposition $\g(D)= \langle D \rangle \oplus \g$,
and its Lie bracket is defined in such a way that
$[D,x]_{\g(D)}=D(x)=-[x,D]_{\g(D)}$, for any $x\in \g$, and
so as to coincide with the Lie bracket $[\,\cdot\,,\,\cdot\,]$ on $\g$
when applied to any pair of elements $x$ and $y$ in $\g$.
This type of extension is actually the semidirect sum,
$\g(D)=\langle D \rangle \ltimes \g$,
and it is called {\bf the extension of $\g$ by the (chosen) derivation} $D\in\Der(\g)$.
In particular, an element $z\in \g(D)$ gets defined in a unique way,
having the property that $[z,x]_{\g(D)}=D(x)$, for any $x\in \g$;
i.e., one may assert that there is a $z\in\g(D)$, such that,
$\ad(z)\vert_{\g} =D$.

\medskip
Now, the two types of extensions just described can be performed
one after the other. Say, start with a Lie algebra $(\g,[\,\cdot\,,\,\cdot\,])$
and a $2$-cocycle $\theta\in\wedge^2\g^*$, so as to produce the
central extension $(\g_\theta,[\,\cdot\,,\,\cdot\,]_\theta)$. Then adjoin
a derivation $D\in\Der(\g_\theta)$ to produce the Lie algebra
$\left(\g(D;\theta),[\,\cdot\,,\,\cdot\,]_{\g(D,\theta)}\right)$,
by further extending the Lie bracket $[\,\cdot\,,\,\cdot\,]_{\theta}$, to the vector space,
$\g(D,\theta):= \langle D \rangle \oplus \g_{\theta}$, letting
$[\,\cdot\,,\,\cdot\,]_{\g(D,\theta)}:\g(D,\theta)\times\g(D,\theta)\to\g(D,\theta)$,
be skew-symmetric and given by, 
$[u,v]_{\g(D,\theta)}=[u,v]_{\theta}$ for any $u$ and $v$ in $\g_{\theta}$,
and $[D,v]_{\g(D,\theta)}=D(v)$ for any $v$ in $\g_{\theta}$.
The process of going from $\g$ to $\g(D,\theta)$ is known as
{\bf the double extension of $\g$ by the pair $(D,\theta)$},
and actually, $\g(D,\theta)=\langle D \rangle \ltimes \g_{\theta}$.
We remark that  the double extension process was first
introduced by V. Ka\c{c} in \cite{Kac}.

\medskip
We shall be dealing with {\it Frobenius Lie algebras\/} 
and with {\it contact Lie algebras\/}. 
Let $(\g,[\,\cdot\,,\,\cdot\,])$ be a Lie algebra.
Given a linear form 
$\varphi\in\g^*$, its associated
{\bf Kirillov form} $B_{\varphi}$ 
is the bilinear form
$B_{\varphi}:\g\times\g\to\F$, defined by,
$$
B_{\varphi}(x,y)=\varphi([x,y]),\quad\text{for all\ }\,x,y\in\g.
$$
Clearly, $B_{\varphi}$ is skew-symmetric.
A Lie algebra $(\g,[\,\cdot\,,\,\cdot\,])$
{\bf is a Frobenius Lie algebra}
if $B_{\varphi}$ is non-degenerate (see \cite{O}).
Equivalently, if $(d\varphi)^n\ne 0$.
We shall also use the notation $(\g,\varphi)$ for
a Frobenius Lie algebra, assuming that $\varphi\in\g^*$
has this property.
In particular, a Frobenius Lie algebra $(\g,\varphi)$
satisfies $\dim(\g)=2n$ for some $n\in\N$,
and it is also {\bf an exact symplectic Lie algebra}
whose symplectic form is $\omega :=B_{\varphi} = -d\varphi$.
Moreover, there is a unique $x_P\in\g$, 
called {\bf the principal element associated to} $\varphi$,
satisfying, $\varphi\,\circ\,\ad(x_P)=\varphi$ (see \cite{O}).
It has been proved in \cite{B-R-S} 
that {\it there is a codimension $1$
unimodular contact ideal\/} $\h$ in $\g$, such that,
$$
\g = \langle x_P\rangle \ltimes \h.
$$
Recall that $(\g,[\,\cdot\,,\,\cdot\,])$ {\bf is a contact Lie algebra},
if $\dim(\g)=2n+1$ for some $n\in\N$,
and if it comes equipped with a
linear form $\alpha\in\g^*$, called {\bf the contact $1$-form},
satisfying $\alpha\wedge(d\alpha)^n\ne 0$. 
Equivalently, in terms of the Kirillov form, there is
an element $\xi\in\g$,
called {\bf the Reeb vector}, satisfying  $\alpha(\xi)=1$,
and $\Rad(B_{\alpha})=\langle \xi \rangle$, where
$\Rad(B_{\alpha})=\{x\in\g\mid B_{\alpha}(x,y)=0,\ \,\forall\ y\in\g\}$
is {\it the radical of\/} $B_{\alpha}$.
We shall also use the notation $(\g,\alpha)$ for
a contact Lie algebra, assuming that $\alpha\in\g^*$
has the required property property, and refer to
$\alpha$ as its {\bf contact structure}.

\medskip
Now, there are, at least,
two different ways of producing
a contact Lie algebra from a Frobenius Lie algebra $(\g,\varphi)$:

\medskip
{\bf (1)} 
Use the symplectic $2$-form $\omega=-d\varphi$
as the defining 2-cocycle for the central extension of $\g$
by a central element $z\notin\g$. 
In this case {\it the contact 1-form\/}
$\alpha\in(\g\oplus\langle z\rangle)^*$,
is $\alpha=z^*$.
\newline
{\bf Note:} For a given $z$ in a Lie algebra $\G$, we follow the convention of 
writing, $z^*\in\G^*$ for a linear functional such that $z^*(z)=1$, and use
later the freedom of choosing the subspace $\Ker(z^*)\subset\G$ as needed.

\medskip
{\bf (2)}
Find a derivation $D\in\Der(\g)$, satisfying $\varphi\,\circ\,D=\varphi$.
Then, $\g(D):=\langle D\rangle \ltimes \g$ becomes a contact Lie algebra
with contact 1-form $\alpha\in\g(D)^*$, given by,
$\alpha=\varphi+\lambda\,D^*$,
for almost any non-zero scalar $\lambda$.
Moreover, in this case there is 
a codimension $1$ ideal $\h$ in $\g$, and $\g$ itself
is a codimension $1$ ideal in $\g(D)$.
Thus, $\g(D) \rhd \g \rhd \h$ (see \cite{B-R-S}).

\medskip
{\bf K\"ahler structures.} In what follows, we shall assume that all
the Lie algebras $\g$ we deal with, are defined over the real field $\R$. A
{\bf K\"ahler structure in a Lie algebra}\/ $(\g,[\,\cdot\,,\,\cdot\,])$ is a triple, 
$(g, J,\omega)$, where,

\medskip
\begin{itemize}
\item [\text{(1)}] $g:\g\times\g\to\R$ is a positive-definite
symmetric bilinear form on $\g$.
\medskip
\item [\text{(2)}] $J:\g\to\g$ is a linear endomorphism satisfying,
\medskip
{
\begin{itemize}
\item [\text{(2.1)}] 
$J\circ J=-\Id_\g$;
\medskip
\item [\text{(2.2)}] 
$g(Jx,y)+g(x,Jy)=0$, for all $x$ and $y$ in $\g$;
\medskip
\item [\text{(2.3)}] 
$N_J(x,y)=0$, for all $x$ and $y$ in $\g$, where 
$N_J$ is the Nijenhuis tensor of $J$, given by,
$$
N_J(x,y)= - [x,y]+[J x,J y]-J[x,J y]-J[J x,y];
$$
\end{itemize}
}
\medskip
\item [\text{(3)}] $\omega:\g\times\g\to\R$, defined through,
$\omega(x,y)=g(x,Jy)$, for all $x$ and $y$ in $\g$, satisfies,
$d\omega=0$.
\end{itemize}

\medskip
{\bf Comments.}
{\bf (1)} In view of (3), it
suffices to define a K\"ahler structure on a given
Lie algebra $(\g,[\,\dot\,,\,\cdot\,])$ in terms of
either, the pair $(g,J)$, or the pair $(\omega,J)$. 
Also, if defined as a pair $(g,J)$, the condition (2.2),
is equivalent to $g(Jx,Jy)=g(x,y)$, for all $x$ and $y$ in $\g$,
given the hypothesis (2.1).
In this context, $g$ is called a 
{\bf Hermitian metric.}
\newline
{\bf (2)} Recall that a linear endomorphism $J\in\End(\g)$
satisfying (2.1) is called an {\bf almost complex structure on} $\g$.
An almost complex structure $J\in\End(\g)$ satisfying (2.3) is called
{\bf a complex structure.}
Finally, a pair $(g,J)$ satisfying (1), (2.1), (2.2) and (3) is called an
{\bf almost K\"ahler structure}.

\medskip
{\bf Sasakian structures.}
A {\bf  Sasakian structure in a Lie algebra\/} $(\g,[\,\cdot\,,\,\cdot\,])$ is a 4-tuple, 
$(g,\xi,\alpha,\Phi)$, where,

\medskip
\begin{itemize}
\item [\text{(1)}] $g:\g\times\g\to\R$ is a positive-definite
symmetric bilinear form on $\g$.
\medskip
\item [\text{(2)}] $\xi\in\g$ and $\alpha\in\g^*$ are dual to each other, 
in the sense that $\alpha(\xi)=1$.
\medskip
\item [\text{(3)}] $\Phi:\g\to\g$ is a linear endomorphism satisfying:
\medskip
{
\begin{itemize}
\item [\text{(3.1)}] 
$\Phi\,\circ\, \Phi=-\Id_\g+\alpha\otimes\xi$;
\medskip
\item [\text{(3.2)}] 
$g(\Phi x,\Phi y)=g(x,y)-\alpha(x)\,\alpha(y)$, for all $x$ and $y$ in $\g$;
\medskip
\item [\text{(3.3)}] 
$g(x,\Phi y)=d\alpha(x,y)$, for all $x$ and $y$ in $\g$;
\medskip
\item [\text{(3.4)}] 
$N_\Phi=-d\alpha\otimes\xi$, where,
$$
N_\Phi(x,y)=\Phi^2[x,y]+[\Phi x,\Phi y]-\Phi[x,\Phi y]-\Phi[\Phi x,y],
$$
for all $x$ and $y$ in $\g$.
\end{itemize}
}
\end{itemize}

\medskip
{\bf Comments.} 
{\bf (1)} The positive-definite
symmetric bilinear form $g:\g\times\g\to\R$ of a K\"ahler 
(or of a Sasakian) structure, need not be
$\ad$-invariant; that is, $g$ does not necessarily
satisfy {\it the invariance condition\/} 
$g([x,y],z)=g(x,[y,z])$, for all $x$, $y$ and $z$ in $\g$.
\newline
{\bf (2)} The conditions
$\Phi\,\circ\, \Phi=-\Id_\g+\alpha\otimes\xi$,
and $\alpha(\xi)=1$, imply that $\Phi(\xi)=0$, $\alpha \circ \Phi = 0$,
and $(\Phi\vert_{\Ker(\alpha)})^2=-\Id_\g$. 
Thus, one might as well use the notation $J:=\Phi\vert_{\Ker(\alpha)}$
(see \cite{Blair}).

\medskip
\begin{Remark}
\rm{
It is proved in \cite{Blair} (Chapter 4) that if 
some kind of structure is given on a Lie algebra $(\g,[\,\cdot\,,\,\cdot\,])$
(or more generally on a differentiable manifold $M$),
there is always  a way to produce an adapted positive-definite inner product
on $\g$ (or Riemannian metric on $M$) compatible with the given structure.
Thus, {\bf (1)} if $(\g,\omega)$ is a symplectic Lie algebra,
there is always a way to produce an almost complex structure
$J\in\End(\g)$ and a positive-definite symmetric bilinear form
$g:\g\times\g\to\R$, satisfying $\omega(x,y)=g(x,Jy)$;
 {\bf (2)} if $(\g,\alpha)$ is a contact Lie algebra,
 there is always a way to produce a linear endomorphism
 $\Phi\in\End(\g)$ and a positive-definite symmetric bilinear form
$g:\g\times\g\to\R$, 
satisfying $d\alpha(x,y)=g(x,\Phi y)$;
or, {\bf (3)} if $(\xi,\alpha,\Phi)$ is a given triple on $\g$
with $\xi\in\g$ and $\alpha\in\g^*$ satisfying, $\alpha(\xi)=1$,
and $\Phi\in\End(\g)$, satisfying $\Phi^2=-\Id_{\g}+\alpha\otimes\xi$, 
and $N_{\Phi}=-d\alpha\otimes\xi$, there is always a way to produce a 
positive-definite symmetric bilinear form
$g:\g\times\g\to\R$, satisfying
$g(\Phi x,\Phi y)=g(x,y)-\alpha(x)\alpha(y)$.
In other words, {\it the existence of an adapted positive-definite
symmetric bilinear form $g:\g\times\g\to\R$, is always guaranteed\/.}
Therefore, it is safe to remove $g$ from the defining data
of K\"ahler and Sasakian structures and take the adapted $g$
for granted.
From now on, we shall also say that $(\g,[\,\cdot\,,\,\cdot\,])$
{\bf is a K\"ahler Lie algebra} (respectively, {\bf a Sasakian Lie algebra})
if $(J,\omega)$ is a K\"ahler structure (respectively, if 
$(\xi,\alpha,\Phi)$ is a Sasakian structure) defined on it.
}
\end{Remark}

\medskip
We now quote two important results that relate
K\"ahler and Sasakian Lie algebras:

\medskip
\begin{Theorem}[See \cite{A-F-V}]
Let $(\xi,\alpha,\Phi)$ be a Sasa\-kian structure on a Lie algebra $(\g,[\,\cdot\,,\,\cdot\,])$
with non-trivial one-dimensional center generated by $\xi$, and let $\h=\Ker(\alpha)$.
Let $J =\Phi\vert_{\h}$ and let $\omega$
be obtained through $d\alpha$ from
the $\h$-component of the Lie bracket on $\g$.
Then, $(J,\omega)$ is a K\"ahler structure on $(\h,[\,\cdot\,,\,\cdot\,]\vert_{\h\times\h})$.
\end{Theorem}

\medskip
In other words, this result states that Sasakian Lie algebras
with non-trivial center generated by $\xi$, must be central extensions
of K\"ahler Lie algebras.

\medskip
\begin{Theorem}[See \cite{Smo}] 
Let $(J,\omega)$ be an almost K\"ahler structure
on the Lie algebra $(\h,[\,\cdot\,,\,\cdot\,])$, 
and let $(\alpha,\xi,\Phi)$ be the corresponding 
{\it contact structure\/} on the central extension
$\h_{\omega}=\h\oplus\langle \xi\rangle$
by the $2$-cocycle defined by $\omega$.
Then the Nijenhuis torsion $N_\Phi$ on $\h_{\omega}$
is given in terms of the Nijenhuis tensor $N_J$ on $\h$ as,
$$
N_\Phi(X,Y)=N_J(x,y)-d\alpha(X,Y)\,\xi,
$$
where, $X=x+\lambda\,\xi$ and $Y=y+\mu\,\xi$, 
for any pair of scalars $\lambda$ and $\mu$ and any pair of elements
$x$ and $y$ in $\h$.
\label{Thm:Kahler-Nijenhuis-Torsion}
\end{Theorem}

\medskip
One concludes from this result that $N_J=0$, if and only if, $N_\Phi=-d\alpha\otimes\xi$.
Therefore, if the almost complex structure $J$ in $\h$ is integrable, then
$\g=\h_{\omega}$, as given in Theorem 
\ref{Thm:Kahler-Nijenhuis-Torsion}, is a Sasakian Lie algebra; in other words:

\medskip
\begin{Cor}
The central extension of a K\"ahler Lie algebra defined by the
$2$-cocycle given by its symplectic form $\omega$, is a Sasakian Lie algebra.
Furthermore, any Sasakian Lie algebra with non trivial center
is obtained as such a central extension
of a K\"ahler Lie algebra.
\end{Cor}

\medskip
{\bf Observation.} 
On the other hand, the next result shows that a
K\"ahler structure $(J,\omega)$ can never be obtained 
on a central extension of a Lie algebra $(\g,[\,\cdot\,,\,\cdot\,])$ 
equipped with a Sasakian structure $(\xi,\alpha,\Phi)$ 
in such a way that its symplectic form $\omega$ be an extension of $d\alpha$
and its almost complex structure $J$ be an extension of $\Phi$. 
More concisely:

\medskip
\begin{Prop}
Let $(\xi, \alpha, \Phi)$ be a Sasakian structure on a given Lie algebra $(\g,[\,\cdot\,,\,\cdot\,])$. 
Then, no central extension of $\g$ produces a K\"ahler structure $(J,\omega)$
whose complex structure satisfies $J\vert_{\Ker({\alpha})} =\Phi$,
and whose symplectic form satisfies $\omega\vert_{\Ker({\alpha})} = - d \alpha$.
\label{Prop:KahlerNoCentralExtensionSasakian}
\end{Prop}

\medskip
\begin{proof}
Assume that the central extension of $(\g,[\,\cdot\,,\,\cdot\,])$ by a $2$-cocycle $\theta$
produces a K\"ahler Lie algebra $(\g_\theta,[\,\cdot\,,\,\cdot\,]_\theta)$. 
Let $z$ be the central element added by the central extension;
i.e., $[x,y]_\theta = [x,y]_\g + \theta(x,y) z$, for any $x$ and $y$ in $\g$.
Assume that the complex structure $J:\g_\theta \to \g_\theta$ is defined by,
$J(x) = \Phi(x)$, for any $x \in \Ker(\alpha)$, and $J(\xi) = z$, and $J(z)= -\xi$.
Since $J$ must be integrable, we have,
$$
\begin{aligned}
0  & = N_J(x,y) = N_{\Phi}(x,y) 
+ d \alpha(x,y) \xi
  \\
& \quad + ( \theta(x,y) + \theta(\Phi(x), \Phi(y)) ) z
 + ( \theta(\Phi(x),y) + \theta(x, \Phi(y)) ) \xi.
\end{aligned}
$$
Therefore, for any $x$ and $y\in \Ker(\alpha)$, one has,
$$
\theta(x,y) + \theta(\Phi(x), \Phi(y))  = 0, \quad \text{and}\quad
\theta(\Phi(x),y) + \theta(x, \Phi(y))  = 0.
$$
Since $\Phi^2 = - \Id$ in $\Ker(\alpha)$, it must be true that
$\theta(x,y)=0$ for any $x$ and $y \in \Ker(\alpha)$. Similarly, 
$N_J(x,\xi)=0$ for any $x \in \Ker(\alpha)$ implies that
$\theta(x,\xi)=0$. It follows that $\theta \equiv 0$. 
This implies that $dz^*= 0$.

\medskip
Now, the hypothesis that the symplectic form on $\g_{\theta}$
should be given by the extension of $d\alpha$, means that,
$$
\omega = d \alpha + \xi^* \wedge z^*.
$$

But then, $d \omega = 0$ if and only if $d \xi^*= 0$, which is impossible
(see Theorem 13 in  \cite{G-R}).
\end{proof}

\medskip
\subsection{Double extensions with integrable $J$'s}

\medskip
Suppose $\h$ is a real Lie algebra with a symplectic form
$\omega:\h\times\h\to\R$ which
is not necessarily exact as for Frobenius Lie algebras. 
Suppose that $J$ is a complex structure on $\h$; i.e., $N_J=0$. 
Suppose that $\g=\h(D,\omega)$
is the double extension of $\h$ by the $2$-cocycle $\omega$ and a derivation $D$ of the
central extension with central element $\xi$. Define $\bar{J}:\g\to\g$ as,
$$
\bar{J}(x)=
J(x),
\ \ \forall\,x\in\h;
\quad \bar{J}(\xi)=D;\quad
\bar{J}(D)=-\xi.
$$

\medskip
\begin{Prop}
Under the hypotheses above, $N_{\bar{J}}=0$ if and only if,
$\bar{J}(Dx)=D(Jx)$, for any $x\in\h$.
\label{Prop:NJ=0}
\end{Prop}

\medskip
\begin{proof}
Clearly, for any $x$ and $y$ in $\h$, one has $N_{\bar{J}}(x,y)=0$.
It is also easy to see that $N_{\bar{J}}(\xi,D)=0$. Similarly, 
$N_{\bar{J}}(x,\xi)=0$ if and only if $\bar{J}(Dx)=D(Jx)$.
Finally, $N_{\bar{J}}(x,D)=0$ if and only if $\bar{J}(D(Jx))=D(x)$, for any $x\in\h$.
Using the fact that $J^2=-\Id_\h$ in the last equality, we conclude that
$N_{\bar{J}}(x,D)=0$ if and only if $\bar{J}(D(x))=D(Jx)$.
\end{proof}

\medskip

\section{Double extensions of K\"ahler and Sasakian Lie algebras}

\medskip
We have seen that {\it in most cases\/} the extension processes
on Lie algebras having some kind of geometric structures,
can be performed in such a way so as to get compatible 
geometric structures in the obtained extensions. 
This is the guiding principle that we have followed in previous works
(e.g., \cite{M-S}, \cite{RV-S} or \cite{RV-S-2}) and that we shall
follow throughout this and the next section.

\medskip
For example, it is now known that the double extension process applied to a given
K\"ahler Lie algebra $(\g,[\,\cdot\,,\,\cdot\,])$, produces again a K\"ahler Lie algebra 
(see \cite{D-M} or \cite{D-M-1}). It is for this reason that in this section
we shall focus our interest in the double extension process, but
applied to Sasakian Lie algebras.


\subsection{Sasakian double extension}

Let $\bar{\alpha}$ be a contact structure
on a Lie algebra $(\g,[\,\cdot\,,\,\cdot\,])$
with Reeb vector $\bar{\xi}\in\g$
and let $\g(\theta, D)=\langle D\rangle\ltimes\g_{\theta}$,
be a double extension of $\g$, where $\g_{\theta}$ is a central
extension of $\g$ by a $2$-cocycle $\theta\in\wedge^2\g^*$,
introducing a central element $z\in\g_{\theta}$, with $z\notin\g$,
and $D\in\Der(\g_{\theta})$. Let $\lambda$ be a non-zero scalar, 
and consider, $\alpha = \bar{\alpha} + \lambda z^*\in\g(\theta, D)^*$.
It has been proved in \cite{RV-S} that if
$$
\alpha (D(z)) \neq 0,
$$
then $(\g(\theta, D), \alpha)$ is a contact Lie algebra. 
As a matter of fact, it suffices to consider the case $\lambda=1$;
that is, $\alpha = \bar{\alpha} + z^*$, since $z$ can be scaled as needed.
The Reeb vector would then have to be of the form,
$\xi=u+a\,\bar{\xi}+b\,z$, 
for some pair of scalars $a$ and $b$ in $\R$,
satisfying $a+b=1$,
and some vector $u\in\Ker(\bar{\alpha})$.
Fix for a moment a vector $w=c\,\bar{\xi}+d\,z\in\Ker(\alpha)$,
for some pair of scalars $c$ and $d$, such that $c+d=0$. 
Therefore, the double extension 
has the vector space decomposition
$\g(\theta, D)=\Ker(\bar{\alpha})\oplus
\langle w, D \rangle \oplus \langle \xi \rangle$.

\medskip
Now, given $(\g,[\,\cdot\,,\,\cdot\,])$, with contact $1$-form $\bar{\alpha}$
and Reeb vector $\bar{\xi}$,
assume there is a linear map $\bar{\Phi} \in \End(\g)$ that turns
$(\bar{\alpha}, \bar{\xi},\bar{\Phi})$ into a Sasakian structure on it.
We may extend $\bar{\Phi}$ to $\Phi \in \End (\g(\theta, D))$
by letting,
$$
\begin{aligned}
\Phi(x)  & := \bar{\Phi}(x), \quad\text{if}\quad x\in \Ker(\bar{\alpha}); 
\\
\Phi(w) & := D, \qquad\ \ \Phi(D) := - w; 
\\
\Phi(\xi) & := 0.
\end{aligned}
\leqno{(1)} 
$$
Observe that the equations
$\Phi(\xi) = 0$ and $\Phi(w) = D$, are equivalent to,
$$
\begin{aligned}
a\,\Phi(\bar{\xi}) + b\,\Phi(z) & =  -\Phi(u), \\
c\,\Phi(\bar{\xi}) + d\,\Phi(z) & = D.
\end{aligned}
$$
If $\delta:= ad - bc\ne 0$, then,
$$
\Phi(\bar{\xi}) = -\frac{1}{\delta} \left( b\,D + d\,\Phi(u)\right), 
\quad\text{and,}\quad
\Phi(z) = \frac{1}{\delta} \left(a\,D + c\,\Phi(u)\right).
$$
We may now proceed to state and prove the main result of this section
using the notation and hypotheses just introduced:

\medskip
\begin{Theorem}
Let $(\bar{\xi}, \bar{\alpha}, \bar{\Phi})$ be a Sasakian structure
on a given Lie algebra $(\g,[\,\cdot\,,\,\cdot\,])$, and let
$\alpha$ be a contact structure on a double extension
$\g(\theta, D)$ of $\g$, with $D\in\Der(\g_\theta)$
and let $\xi\in\g(\theta,D)$ be its Reeb vector. 
Let $z$ be the central element introduced by the given $2$-cocycle $\theta$.
Write $\alpha = \bar{\alpha} + z^*$, and assume that the Reeb vector 
has the form, $\xi=u+a\,\bar{\xi}+b\,z$, 
for some pair of scalars $a$ and $b$ in $\R$,
satisfying $a+b=1$,
and some vector $u\in\Ker(\bar{\alpha})$.
Also assume that $w=c\,\bar{\xi}+d\,z\in\Ker(\alpha)$,
so that $c+d=0$.
Define $\Phi\in\End(\g(\theta, D))$ in terms of $\bar{\Phi}$ as in {(1)} above.
Then, $(\g(\theta, D), \xi, \alpha, \Phi)$ is a Sasakian Lie algebra
if and only if,
\begin{enumerate}
\item 
$\theta( \bar{\Phi} x,y) + \theta (y, \bar{\Phi} y) = 0$, for any
$x$ and $y$ in $\Ker(\bar{\alpha})$;
\item $u$ and $\bar{\xi}$ belong to $\Rad(\theta)$;
\item $D\circ \bar{\Phi} = \Phi \circ D$, on $\Ker(\bar{\alpha})$;
\item $\ad_{\g}(u) =  - \bar{\Phi} \circ \ad_{\g}(u) \circ \bar{\Phi}$;
\item $c\, [\bar{\xi}, u] + \Phi(D(\xi)) = 0$.
\end{enumerate}
\end{Theorem}

\medskip
We shall prove this Theorem by a series of observations and
calculations stated as lemmas. We start by noting that
what needs to be proved is that,
$$
\Phi^2 = - \Id + \alpha \otimes \xi,
$$
so that, $N_\Phi=-d\alpha\,\otimes\,\xi$,
if and only if,
$$
M(U,V):= - [U,V] + [\Phi U, \Phi V] - \Phi( [\Phi U, V]) - \Phi([U, \Phi V]) = 0,
$$
for any $U$ and $V$ in $\g(\theta,D)$.
What we shall find in the following lemmas is the set of conditions
(1)-(5) above, in terms of what is needed for
$M(U,V)$ to vanish for any $U$ and $V$ in $\g(\theta,D)$. 
Thus we start with the following:

\medskip
\begin{Lemma}
For any $x$ and $y$ in $\Ker(\bar{\alpha})$, 
$$
M(x,y) = 0
\quad\Longleftrightarrow\quad
\theta( \bar{\Phi} x,y) + \theta (x, \bar{\Phi} y) = 0.
$$
\end{Lemma}

\medskip
\begin{proof}
Take $U=x$ and $V=y$ with $x$ and $y$ in $\Ker(\bar{\alpha})$. Then,
it is straightforward to see that,
$$
\begin{aligned}
M(x,y)  
& =
  - [x,y] + [\Phi x, \Phi y] - \Phi([\Phi x, y]) - \Phi( [x, \Phi y]) \\
&  = 
 - [x,y]_{\g} - \theta(x,y)\,z + [\bar{\Phi} x, \bar{\Phi} y]_\g + \theta( \bar{\Phi} x, \bar{\Phi} y)\,z \\
& 
\quad
 -\bar{\Phi} ( [\bar{\Phi} x, y]_{\g} ) + \theta( \bar{\Phi}x, y) \,\Phi (z) \\
& 
\quad
 -\bar{\Phi} ( [x, \bar{\Phi}y]_{\g} ) + \theta(x, \bar{\Phi} y) \,\Phi (z) \\
& = 
  - [x,y]_\g + [\bar{\Phi} x, \bar{\Phi} y]_\g  
   -\bar{\Phi} ( [\bar{\Phi} x, y]_{\g} ) -\bar{\Phi} ( [x, \bar{\Phi} y]_{\g} )
\\
&  
\quad
+ \left(\,\theta( \bar{\Phi} x, \bar{\Phi} y) - \theta(x,y) \,\right)\,z  
\\
&  
\quad
+ \left(\,\theta( \bar{\Phi}x, y) + \theta(x, \bar{\Phi} y) \,\right)\, \Phi(z).
\end{aligned}
$$
Now, $- [x,y]_\g + [\bar{\Phi} x, \bar{\Phi} y]_\g  
-\bar{\Phi} ( [\bar{\Phi} x, y]_{\g} ) -\bar{\Phi} ( [x, \bar{\Phi} y]_{\g} )=0$,
because $x$ and $y$ lie in $\Ker(\bar{\alpha}) \subset \g$, and $\g$ is
Sasakian. On the other hand, $\theta( \bar{\Phi} x, \bar{\Phi} y) - \theta(x,y)=0$,
if and only if,
$\theta( \bar{\Phi}x, y) + \theta(x, \bar{\Phi} y)=0$.
In particular, it suffices that one or the other vanishes to conclude the 
truth of the statement.

\end{proof}

\medskip
\begin{Lemma}
Assume that $w=c\,\bar{\xi} + d\, z \in \Ker(\alpha)$, for
appropriate choices of the scalars $c$ and $d$, and let
$x$ be in $\Ker(\bar{\alpha})$. Then,
$$
M(x,w)=0
\quad\Longleftrightarrow\quad
\begin{cases}
(a) \ \bar{\xi} \in \Rad(\theta),\ \ \text{and,}
\\
(b) \ D\circ\bar{\Phi} = \Phi \circ D,\,\,\text{on\ }\Ker(\bar{\alpha}).
\end{cases}
$$ 
\end{Lemma}

\medskip
\begin{proof}
Asume that $x \in \Ker(\bar{\alpha})$ and write, $w= c\,\bar{\xi} + d\,z \in \Ker(\alpha)$.
Then,
$$
\begin{aligned}
M(x, c\,\bar{\xi} + d\,z)  
& = 
-[x, c\,\bar{\xi} + d\,z] + [\Phi(x), \Phi(c\,\bar{\xi} + d\,z)]
\\
& \quad
- \Phi( [\Phi(x), c\,\bar{\xi} + d\,z])
- \Phi( [x, \Phi(c\,\bar{\xi} + d\,z)]) 
\\
& =
- c\,[x, \bar{\xi}] + c\,[\Phi(x), \Phi(\bar{\xi} )] + d\,[\Phi(x), \Phi( z)] 
\\
& \quad
 - c\,\Phi( [\Phi(x), \bar{\xi} ]) - c\,\Phi( [x, \Phi(\bar{\xi} )]) - d\,\Phi( [x, \Phi(z)])
 \\
& =
 -c\,\left( [x, \bar{\xi}]_\g + \theta(x, \bar{\xi})\,z\right) 
 + c\,\left( [\bar{\Phi}(x), -\frac{1}{\delta}( b\,D + d\,\Phi(u))]\right) 
 \\
& \quad
 + d\, [\bar{\Phi}(x), \frac{1}{\delta} (a\,D + c\,\Phi(u)) ]
 - c\,\Phi \left( [\bar{\Phi}(x), \bar{\xi}]_\g + \theta( \bar{\Phi}(x), \bar{\xi})\,z\right)  
\\
& \quad
- c\, \Phi\left(  [x, -\frac{1}{\delta}( b\,D + d\,\Phi(u)] \right)
- d\, \Phi \left( [x, \frac{1}{\delta} ( a\,D + c\,\Phi(u) ] \right) 
\\
& = 
-c\, ([x, \bar{\xi}]_\g + \bar{\Phi}( [\bar{\Phi}(x), \bar{\xi}]_\g) 
+ \left( \frac{ad}{\delta} - \frac{bc}{\delta} \right) [\bar{\Phi}(x), D] 
\\
& \quad
 + \left( \frac{bc}{\delta} - \frac{ad}{\delta} \right) \Phi( [x, D]) 
 - \frac{cd}{\delta} [\bar{\Phi}(x), \Phi(u)] 
 \\
& \quad
 + \frac{cd}{\delta}  [\bar{\Phi}(x), \Phi(u)] 
 + \frac{cd}{\delta} \Phi( [x, \Phi(u)] )- \frac{cd}{\delta}  \Phi( [x, \Phi(u)]) 
 \\
& \quad
 - c\,\theta(x, \bar{\xi})\,z - c\,\theta (\bar{\Phi}(x), \bar{\xi} )\,\Phi(z)
 \\
 & = 
 -c\, ([x, \bar{\xi}]_\g + \bar{\Phi}( [\bar{\Phi}(x), \bar{\xi}]_\g) 
+ [\bar{\Phi}(x), D] - \Phi( [x, D])
\\
& \quad
 - c\,\theta(x, \bar{\xi})\,z - c\,\theta (\bar{\Phi}(x), \bar{\xi} )\,\Phi(z).
\end{aligned}
$$
Now, $\bar{\Phi}( [\bar{\Phi}(x), \bar{\xi}]_\g) = c\, [x, \bar{\xi}]_\g$
because $\g$ is Sasakian. Whence, 
$$
M(x,c\,\bar{\xi}+d\,z)=
[\bar{\Phi}(x), D] - \Phi( [x, D])
- c\,\theta(x, \bar{\xi})\,z - c\,\theta (\bar{\Phi}(x), \bar{\xi} )\,\Phi(z).
$$
Therefore, if $\bar{\xi} \in \Rad(\theta)$, 
we conclude that $M(x, c \bar{\xi} + d z)  = 0$, if and only if, 
$D(\bar{\Phi} (x)) = \Phi (D (x))$, for all $x\in \Ker(\bar{\alpha})$.
\end{proof}

\medskip
\begin{Lemma}
For any $x\in \Ker(\bar{\alpha})$, 
$$
M(x,D)=0
\quad\Longleftrightarrow\quad
\begin{cases}
(a) \ \bar{\xi} \in \Rad(\theta),\ \ \text{and,}
\\
(b) \ D\circ \bar{\Phi} = \Phi \circ D,\,\, \text{on\ } \Ker(\bar{\alpha}).
\end{cases}
$$
\end{Lemma}

\medskip
\begin{proof}
Assume $x\in \Ker(\bar{\alpha})$, and write $\Phi(D) = - c\,\bar{\xi} - d\,z$, as before.
It follows that,
$$
\begin{aligned}
M(x,D) 
& = 
 - [x, D] + [\Phi(x), \Phi(D)] - \Phi ([\Phi(x), D]) - \Phi([x, \Phi(D)] 
 \\
& =
  D(x) - c\,[\bar{\Phi}(x), \bar{\xi}] - \Phi([\bar{\Phi}(x), D]) + c\,\Phi( [x, \bar{\xi}]) 
  \\
& =
   D(x) + \Phi([\bar{\Phi}(x), D])  -  c\,
   \left( [\bar{\Phi}(x), \bar{\xi}]_\g + \theta( \bar{\Phi}(x), \bar{\xi}) \,z \right) 
\\
& \quad
+  c\,\Phi ( [x, \bar{\xi}]_\g + \theta( x, \bar{\xi})\,z )
\\
& = 
 D(x) - \Phi ([\Phi(x), D]) 
 \\
& \quad
 + c\,( \Phi ( [x, \bar{\xi}]_\g ) -  [\bar{\Phi}(x), \bar{\xi}] ) 
 \\
& \quad
 + c\, (  \theta( x, \bar{\xi})\,\Phi(z) -  \theta( \bar{\Phi}(x), \bar{\xi})\,z.
\end{aligned}
$$
Since $\g$ is Sasakian, $\Phi ( [x, \bar{\xi}]_\g ) =  [\bar{\Phi}(x), \bar{\xi}]$,
and we are left with,
$$
M(x,D) = D(x) - \Phi ([\Phi(x), D])
+ c\, (  \theta( x, \bar{\xi})\,\Phi(z) -  \theta( \bar{\Phi}(x), \bar{\xi})\,z.
$$
It follows as before that if $\bar{\xi}\in\Rad(\theta)$, then
$M(x,D)= 0$, if and only if, $D(x) = \Phi ([\Phi(x), D])$; equivalently,
if and only if, $\Phi(D(x)) = D(\bar{\Phi}(x))$.
\end{proof}

\medskip
\begin{Lemma}
Let $x$ be in $\Ker(\bar{\alpha})$, and write $\xi=u + a\,\bar{\xi} + b\,z$,
with $u\in\Ker(\bar{\alpha})$, as before. Then,
$$
M(x,\xi)=0
\quad\Longleftrightarrow\quad
\begin{cases}
(a) \ u\text{\,\,and\,\,\,}\bar{\xi} \in \Rad(\theta),\ \ \text{and,}
\\
(b) \ \ad_{\g}(u) =  - \bar{\Phi} \circ \ad_{\g}(u) \circ \bar{\Phi}.
\end{cases}
$$ 
\end{Lemma}

\medskip
\begin{proof}
Assume $x\in \Ker(\bar{\alpha})$. We know in that, $\Phi(\xi)=0$. 
Whence,
$$
\begin{aligned}
M(x, \xi) 
& =
 - [x, \xi] + [\Phi(x), \Phi(\xi)] - \Phi([\Phi(x), \xi]) - \Phi( [x, \Phi(\xi)]) 
 \\
& =
 -[x,u]  - a\,[x, \bar{\xi}] - \Phi( [\bar{\Phi}(x), u + a\,\bar{\xi} + b\,z])  
 \\
& = 
 - ([x,u]_\g + \theta(x,u)\,z) - a\,([x,\bar{\xi}]_\g + \theta(x, \bar{\xi})\,z) 
 \\
& \quad
 -  \Phi( [\bar{\Phi}(x), u]) - a\,\Phi( [\bar{\Phi}(x), \bar{\xi}])) 
 \\
& =
 - ([x,u]_\g + \theta(x,u)\,z) - a\,([x,\bar{\xi}]_\g + \theta(x, \bar{\xi})\,z) 
 \\
& \quad
 - \Phi( [\bar{\Phi}(x), u]_\g + \theta(\bar{\Phi}(x), u)\,z)  - a\,\Phi( [\bar{\Phi}(x), \bar{\xi}])) 
 \\
& =
 - [x,u]_\g - \Phi( [\bar{\Phi}(x), u]_\g )
 \\
& \quad
 - (\theta(x,u)\,z  + \theta(\bar{\Phi}(x), u)\,\Phi(z)) - a\,\theta(x, \bar{\xi})\,z 
 \\
& \quad
 - a\,([x,\bar{\xi}]_\g + \Phi( [\bar{\Phi}(x), \bar{\xi}]) ).
\end{aligned}
$$
Since $\g$ is Sasakian, $[x,\bar{\xi}]_\g =-\Phi( [\bar{\Phi}(x), \bar{\xi}])$.
Under the assumption that $u$ and $\bar{\xi}$ belong to $\Ker(\theta)$,
it follows that 
$M(x, \xi)=0$, if and only if,
$- [x,u]_\g - \Phi( [\bar{\Phi}(x), u]_\g )=0$.
\end{proof}

\medskip
\begin{Lemma}
If $w= c\,\bar{\xi}+ d\, z \in \Ker(\alpha)$, then $M(w,D)=0$.
\end{Lemma}

\medskip
\begin{proof}
Since $\Phi(w)=D$, and $\Phi(D)=-w$, it trivially follows that $M(w,D)=0$.

\end{proof}

\medskip
\begin{Lemma}
Given $w = c\,\bar{\xi}+ d\,z \in \Ker(\alpha)$, and $\xi=u + a\,\bar{\xi} + b\,z$,
with $u\in\Ker(\bar{\alpha})$,
$$
M(w,\xi)=0
\quad\Longleftrightarrow\quad
\begin{cases}
(a) \ \bar{\xi} \in \Rad(\theta),\ \ \text{and,}
\\
(b) \ c\,[\bar{\xi}, u]_\g + \Phi(D(\xi)) = 0.
\end{cases}
$$ 
\end{Lemma}

\medskip
\begin{proof}
Using the fact that $\Phi(\xi)=0$, and the fact that $\Phi(w) = D$, we have,
$$
\begin{aligned}
M(w, \xi) 
& = 
- [w, \xi] + [\Phi(w), \Phi(\xi)] - \Phi( [\Phi(w), \xi]) - \Phi([w, \Phi(\xi)])
\\
& =
- [c\,\bar{\xi} + d\,z, u + a\,\bar{\xi} + b\,z]  -\Phi( [D, \xi]) 
\\
& = 
- c\,[\bar{\xi}, u] - \Phi( D(\xi) ) 
\\
& =
- c\,( [\bar{\xi}, u]_\g + \theta(\bar{\xi}, u)\,z ) - \Phi(D(\xi)).
\end{aligned}
$$
If $\bar{\xi}\in \Ker(\theta)$, the assertion follows.

\end{proof}

\medskip
\begin{Lemma}
Write,
$\xi=u + a\,\bar{\xi} + b\,z$,
with $u\in\Ker(\bar{\alpha})$. Then,
$$
M(D,\xi)=0
\quad\Longleftrightarrow\quad
\begin{cases}
(a) \ \bar{\xi} \in \Rad(\theta),\ \ \text{and,}
\\
(b) \ c\,[\bar{\xi}, u] + \Phi(D(\xi)) = 0.
\end{cases}
$$ 
\end{Lemma}

\medskip
\begin{proof}
Using the fact that, $\Phi(D)= -w$, and $\Phi(\xi)=0$, we have,
$$
\begin{aligned}
M(D,\xi) 
& =
 - [D, \xi] + [\Phi(D), \Phi(\xi)] - \Phi( [\Phi(D), \xi] ) - \Phi( [ D, \Phi(\xi)] ) 
 \\
& =
 - D(\xi) + \Phi( [ c\,\bar{\xi} + d\,z, u + a\,\bar{\xi} + b\,z]) 
 \\
& =
 - D(\xi) + c\,\Phi( [\bar{\xi}, u] )  
 \\
& =
 - D(\xi) + c\,\Phi( [\bar{\xi}, u]_\g + \theta(\bar{\xi}, u)\,z ).
\end{aligned}
$$
Therefore, for $\bar{\xi}\in \Rad(\theta)$, one has,
$M(D,\xi)=0$, if and only if, $- D(\xi) + c \Phi( [\bar{\xi}, u]_\g )=0$, 
and this equation is equivalent to, $c\,[\bar{\xi}, u] + \Phi(D(\xi)) = 0$.

\end{proof}

\medskip

We have thus proved {\bf Theorem 3.1}. {\bf qed}


\medskip

\section{Extensions of Frobenius-K\"ahler and Sasakian Lie algebras}

We are now interested in the following problem:
Given a K\"ahler (respectively, Sasakian) structure 
on a Lie algebra $(\g,[\,\cdot\,,\,\cdot\,])$,
under what conditions can we 
conclude that the extension
$\g(D)=\langle D\rangle\oplus\g$, with $D\in\Der(\g)$
carries a Sasakian (respectively, K\"ahler) structure?
This section gives precisely these conditions. 


\medskip
\subsection{From Frobenius-K\"ahler to Sasakian Lie algebras} 

Let $(J,\omega)$ be a K\"ahler structure on the
Lie algebra $(\g,[\,\cdot\,,\,\cdot\,])$. 
One says that, $\g$ {\bf is a Frobenius-K\"ahler Lie algebra}
if there is some $\alpha\in\g^*$, such that $\omega=-d\alpha$,
thus turning the pair $(\g,\alpha)$ into a Frobenius Lie algebra.

\medskip
Now, let $D\in\Der(\g)$ be such that $\alpha\circ D = 0$ and $[D,J]=0$, and let
$\g(D) := \langle D\rangle \ltimes \g$,
be the extension of $\g$ by $D$. Let $\xi\in\g(D)$ be 
such that $D(x)=[\xi,x]$, for any $x\in\g$.

\medskip

With the notation and hypotheses just introduced, we shall
now state the following:

\medskip
\begin{Theorem} \label{FKtoSasaki}
Let $(J,\omega)$ be a Frobenius-K\"ahler structure on the
Lie algebra $(\g,[\,\cdot\,,\,\cdot\,])$, 
with $\omega=-d\alpha$ for some $\alpha\in\g^*$.
Let $D\in \Der(\g)$ be chosen so that
$\alpha \circ D = 0$ and $D\circ J = J \circ D$. 
Then, the triple $(\xi,\bar{\alpha},\Phi)$ turns $\g(D)$
into a Sasakian Lie algebra.
\end{Theorem}

\medskip
\begin{proof}
Let $\bar{\alpha}:= \alpha + \xi^*\in\g(D)^*$. Then,
 $(\g(D), \bar{\alpha})$ is a contact Lie algebra.
This follows easily because $d\xi^*=0$ and $-d\alpha=\omega$.
Therefore, $\bar{\alpha} \wedge (d \bar{\alpha})^n  = \xi^* \wedge \omega^n \neq 0$.
Now define $\Phi : \g(D) \to \g(D)$ by means of
$\Phi(x) = J(x) + \bar{\alpha}(x) \xi$, for any $x \in \g$ and $\Phi(\xi)=0$.
This definition of $\Phi$ turns $\xi$ into the Reeb vector
associated to $\bar{\alpha}$. Now, $J^2=-\Id$, implies that,
$\Phi^2 = -\Id + \bar{\alpha} \otimes \xi$.
Thus, all what we need to prove is that,
$$
N_{\Phi} = - d \bar{\alpha} \otimes \xi.
$$
To prove this fact, recall first that
for any pair $u$ and $v$ in $\g(D)$,
$$
N_{\Phi}(u,v) = \Phi^2 [u,v] + [\Phi(u), \Phi(v)] - \Phi([\Phi u, v]) - \Phi([u, \Phi v]).
$$
Since, $\Phi^2 [u,v] = - [u,v] + \bar{\alpha}([u,v]) \xi$, the proof of (2)
comes down to show that, for any $u$ and $v$ in $\g(D)$,
$$
- [u,v] + [\Phi(u), \Phi(v)] - \Phi([\Phi u, v]) - \Phi([u, \Phi v]) = 0.
$$
For $u=x$, and $v=y$, both in $\g$, one has,
$$
\begin{aligned}
- [x,y]  & + [\Phi(x), \Phi(y)] - \Phi([\Phi x, y]) - \Phi([x, \Phi y]) = \\
           & = [Jx, Jy] + \bar{\alpha}(x) [\xi, Jy] + \bar{\alpha}(y) [Jx, \xi] \\
           & \quad - J([Jx,y]) - \bar{\alpha}([Jx,y]) \xi - \bar{\alpha}(x) J([\xi, y]) - 
                \bar{\alpha}(x) \bar{\alpha}([\xi, y]) \xi \\
           & \quad - J([x,Jy]) - \bar{\alpha}([x,Jy]) \xi - \bar{\alpha}(y) J([x, \xi]) -
                \bar{\alpha}(y) \bar{\alpha}([x, \xi]) \xi \\
       & = N_J(x,y) + \bar{\alpha}(x) [D,J](y) - \bar{\alpha}(y) [J,D](x) \\
       & \quad + \left(d\alpha (Jx,y) + d\alpha (x, Jy) + \bar{\alpha}(x) d\bar{\alpha}(\xi, y)
       + \bar{\alpha}(x) d\bar{\alpha}(x, \xi) \right) \xi \\
       & = 0.
\end{aligned}
$$
On the other hand, for $u=x\in\g$, and $v=\xi$, one has,
$$
\begin{aligned}
N_{\Phi}(x,\xi)  & = \Phi^2 [x,\xi] + [\Phi(x), \Phi(\xi)] - \Phi([\Phi x, \xi]) - \Phi([x, \Phi \xi]) \\
                        & = - [x, \xi] - \Phi( [Jx  +\bar{\alpha}(x) \xi, \xi] ) \\
                        & = - [x, \xi] - \Phi ([Jx, \xi]) = - [x, \xi] - J( [Jx,\xi]) \\
                        & = - [x, \xi] + J(D(Jx)) = - [x, \xi] -[\xi, x] \\
                        & = 0.
\end{aligned}
$$
\end{proof}


\subsection{From Sasakian to  Frobenius-K\"ahler Lie algebras}

\medskip
\begin{Theorem}
Let  $(\xi,\alpha,\Phi)$ be a Sasakian structure
on the Lie algebra $(\g,[\,\cdot\,,\,\cdot\,])$.
Let $D\in\Der(\g)$ be such that $\alpha \circ D = \alpha$, and 
assume that $[\Phi, D]\vert_{\Ker(\alpha)}=0$. Consider the Lie algebra
$\g(D)=\langle D \rangle \ltimes \g$, and let $x_P$ be
the unique element in $\g(D)$ such that, 
$[x_P, x] = D(x)$, for any $x\in \g$. Now define
$J: \g(D) \to \g(D)$ by letting 
$J(x)= \Phi(x) + \alpha(x) x_P$, and $J(x_P) = - \xi$.
Further define, $\bar{\alpha} = \alpha + x_P^* \in \g(D)^*$.
Then, $(\g(D),\bar{\alpha})$
is a Frobenius-K\"ahler Lie algebra.
\end{Theorem} 

\medskip
\begin{proof}
The condition $\alpha \circ D = \alpha$ implies that 
$\bar{\alpha} = \alpha + x_P^*$ is a Frobenius structure 
on $\g(D)$ (see \cite{B-R-S}), 
making $\omega = -d\bar{\alpha}$ symplectic.
We shall now prove that for the defined $J: \g(D) \to \g(D)$,
$N_J(u,v) = 0$, for any $u$ and $v$ in $\g(D)$. 
Observe that $\g(D)$ has the vector space decomposition,
$$
\g(D)=
\Ker(\alpha) \oplus \langle \xi \rangle \oplus \langle x_P \rangle.
$$
Assume first that $u=x$ and $v=y$ lie in $\Ker(\alpha)$. Then,
$$
\begin{aligned}
N_J(x,y) & = - [x,y]   + [J x, J y] - J([J x, y]) - J ([x, J y])  \\
           & = \Phi^2 ([x,y]) - \alpha([x,y]) \xi + [\Phi(x) + \alpha(x) x_P, \Phi(y) + \alpha(y) x_P] \\         
          & \quad - J([ \Phi(x) + \alpha(x) x_P, y]) - J([x, \Phi(y) + \alpha(y) x_P]) \\
          & = \Phi^2 ([x,y]) - \alpha([x,y]) \xi + [\Phi(x), \Phi(y)] + \alpha(y) [\Phi(x), x_P] \\
          & \quad + \alpha(x) [x_P, \phi(y)]  - \Phi( [\Phi(x), y]) - \alpha([\Phi(x),y]) x_P \\
          & \quad - \alpha(x) \Phi( [x_P, y])  - \alpha(x) \alpha([x_P, y]) x_P - \Phi([x, \Phi(y)]) \\
          & \quad - \alpha([x, \Phi(y)]) x_P - \alpha(y) \Phi( [x, x_P]) - \alpha(y) \alpha([x,x_P]) x_P \\
          & = N_\Phi(x,y) + d\alpha (x,y)  + \alpha(y) ( [\Phi(x), x_P] - \Phi( [x, x_P]))  \\
          & \quad + \alpha(x) ( [x_P, \Phi(y)] - \Phi( [x_P, y]) ) \\
          & \quad - ( \alpha(y) \alpha([x,x_P]) + \alpha(x) \alpha([x_P,y])) x_P \\
          & \quad - ( \alpha( [x, \Phi(y)] ) + \alpha( [\Phi(x), y]) x_P \\
       & = 0.
\end{aligned}
$$
On the other hand, now assume 
$u=x\in \Ker(\alpha)$ and $v=x_P$. Then,
$$
\begin{aligned}
N_J(x,x_P) & = - [x,x_P]   + [J x, J x_P] - J([J x, x_P]) - J ([x, J x_P])  \\
 & = [x,x_P] + [\Phi(x) + \alpha(x) x_P, -\xi] - J( [ \Phi(x) + \alpha(x) x_P, x_P])        \\
 & \quad - J ([x, - \xi]) \\
 & = - [x, x_P] - [\Phi(x), \xi] - \alpha(x) \xi  - \Phi([\Phi(x), x_P]) - \alpha( [\Phi(x), x_P]) x_P \\
 & \quad + \Phi([x, \xi]) + \alpha ([x, \xi]) x_P \\
 & = - ( [x,x_P] + \Phi([\Phi(x), x_P]) ) - ( [\Phi(x), \xi] - \Phi( [x, \xi] )) \\
 & \quad + ( \alpha( [x, \xi] ) - \alpha ([\Phi(x), x_P]) ) x_P - \alpha(x) \xi  \\
 & = 0.
\end{aligned}
$$
Finally, since $[x_P, \xi] = \xi$, $J(\xi) = x_P$, and $J(x_P)= -\xi$, 
one trivially has, $N_J(x_P, \xi) = 0$.
\end{proof}

\medskip
\begin{Cor}
Let $(\alpha, J)$ be a Frobenius-K\"ahler structure on the Lie algebra
$(\g,[\,\cdot\,,\,\cdot\,])$ with symplectic form  $\omega=-d\alpha$.
Let $\g\oplus\langle\xi\rangle$
be the central extension of $\g$ by a central element $\xi$ via the $2$-cocycle
defined by $\omega$. Suppose there is a derivation 
$D\in\Der(\g\oplus\langle\xi\rangle)$ satisfying $[D,J]=0$,
and $\xi^*\circ D=\xi^*$.
Then the double extension $\g(D,\omega)$ of $\g$
is a Frobenius- K\"ahler Lie algebra.
\end{Cor}

\medskip

{\bf Notation.}
Let $(\alpha, J)$ be a Frobenius-K\"ahler structure 
on a Lie algebra $(\g,[\,\cdot\,,\,\cdot\,])$
and let $x_P$ be its principal vector.
We know that $\g = \langle x_P \rangle \ltimes \h$,
where $\h$ is a codimension $1$ unimodular contact ideal
(see \cite{M-S}) whose contact structure is just the
restriction $\alpha\vert_{\h}$. 
Let $i: \h\hookrightarrow\g$ be the natural inclusion and
let $\xi$ be the Reeb vector associated to $\bar{\alpha}:=i^*\alpha$.
Then $\h$ can be decomposed as $\h=\langle\xi\rangle\oplus\Ker(\bar{\alpha})$
and one may define the linear endomorphism $\Phi:\h\to\h$, by letting
$\Phi(\xi)=0$ and $\Phi(x)=J(x)$ for any $x\in \Ker(\bar{\alpha})$.
Arguing as before, we can state the following:

\medskip
\begin{Theorem}
With the hypotheses of the last paragraph, 
$(\xi,\alpha,\Phi)$ is a Sasakian structure on $\h$
if and only if $[\ad(\xi),\Phi]=0$.
\end{Theorem}

\medskip
\begin{Remark}
\rm{
In any Frobenius-K\"ahler Lie algebra it is true that
$[x_P, \xi]=\xi$.
If $J(\xi) = x_P$, it follows that $[\ad(\xi), \Phi]=0$ if and only if
$[\ad(x_P), \Phi]=0$, but $[\ad(x_P), \Phi]$ is zero only 
on $\Ker(\bar{\alpha})$.
}
\end{Remark}
\medskip
\begin{Cor} With the notation and hypotheses above,
$\h$ is a Sasakian Lie algebra if and only if 
$[\ad(x_P), \Phi]=0$ on $\Ker(\bar{\alpha})$.
\end{Cor}


\section{Examples}

Here are two examples of what we have done so far.

\medskip
Let $\h_3=\operatorname{Span}_{\R}\{e_1,e_2,e_3\}$ be the
Heisenberg Lie algebra with its Lie bracket defined by
$[e_1, e_2]= e_3$,
and $[e_i, e_3]= 0$, for $i=1$ or $2$.
Let $\{e^1,e^2,e^3\}$ be the dual basis to $\{e_1,e_2,e_3\}$.
It is well known that $\h_3$ is a Sasakian Lie algebra.
Its contact structure is defined by $\alpha = e^3 \in \h_3^*$.
Its Reeb vector is defined by $\xi=e_3$, and the linear
endomorphism $\Phi: \h_3 \to \h_3$ may be defined by 
$\Phi(e_1)=e_2$, $\Phi(e_2)=-e_1$, and $\Phi(e_3)=0$.

\medskip
Let $D\in \Der(\h_3)$ be defined
in such a way that $\alpha \circ D = \alpha$, for $\alpha=e^3$.
Write, $D(e_j)=D_{1j}e_1+D_{2j}e_2+D_{3j}e_3$,
for $1\le j\le 3$. Then,
$$
\alpha \circ D = \alpha
\quad\Longleftrightarrow\quad
D_{3j}=\delta_{3j}.
$$
Moreover, since $D\in \Der(\h_3)$, one may compute
$D([e_i,e_j])$ for $1\le i<j\le 3$ and conclude that
$D_{11}+D_{22}=1$, and $D_{13}=0=D_{23}$. 
In particular, we may fix $D\in \Der(\h_3)$ satisfying
$\alpha \circ D = \alpha$ by letting,
$$
D(e_1)=\displaystyle{\frac{1}{2} } \,e_1,\qquad
D(e_2)=\displaystyle{\frac{1}{2}}\,e_2,\qquad
D(e_3)=\,e_3.
$$
Let $e_4$ be the
unique element in $\h_3(D)$ satisfying $[e_4, e_i]= D(e_i)$ for
any $1\le i\le 3$, so that $x_P=e_4$. Whence, $[x_P,\xi]=\xi$;
equivalently, $[e_4,e_3]=e_3$.

\medskip
Now, compute $d\alpha$. It is clear that,
$$
\aligned
de^3(e_1,e_2)&=-e^3([e_1,e_2]) = -1,
\\
de^3(e_3,e_4)&=-e^3([e_3,e_4]) = 1,
\\
de^3(e_i,e_j)&=0,\quad\text{for all others\ \ }1\le i<j\le 4.
\endaligned
$$
Therefore,
$$
d\alpha=de^3 = e^1\wedge e^2 - e^3\wedge e^4.
$$
It follows that $\h_3(D)$ is a Frobenius-K\"ahler Lie algebra with
symplectic (exact) form, 
$\omega = -d\alpha = -e^1\wedge e^2 + e^3\wedge e^4$.

\medskip
{\bf Remark.}
We are following the convention that, for any Lie algebra $\g$,
each $\chi\in\g^*$ defines a degree $-1$ derivation $i_\chi\in\Der(\wedge\g)$,
with the property that, $i_\chi(x)=\chi(x)$ for any $x\in\g$. Then,
using the universal property of the exterior algebra $\wedge\g^*$,
one may identify $i_{\chi_1}\circ\cdots\circ i_{\chi_k}$ with 
${\chi_1}\wedge\cdots\wedge{\chi_k}\in\wedge\g^*$, for a collection
of elements $\chi_j\in\g^*$, $1\le j\le \dim(\g)$. In particular,
$$
\aligned
(\chi_1\wedge \chi_2) (x_1\wedge x_2) 
& = 
(i_{\chi_1}\circ i_{\chi_2}) (x_1\wedge x_2) 
\\
& =
i_{\chi_1}(i_{\chi_2}(x_1)\,x_2 - i_{\chi_2}(x_2)\,x_1)
\\
& =
i_{\chi_1}({\chi_2}(x_1)\,x_2 - {\chi_2}(x_2)\,x_1)
\\
& =
\chi_2(x_1)\,\chi_1(x_2) - \chi_2(x_2)\,\chi_1(x_1) 
\\
& =
-\det
{\begin{pmatrix}
\chi_1(x_1) & \chi_1(x_2) \\
\chi_2(x_1) & \chi_2(x_2) \\
\end{pmatrix}}.
\endaligned
$$

\medskip
We now want to define a complex structure $J\in\End(\h_3(D))$
by extending $\Phi(e_1)=J(e_1)=e_2$, and $\Phi(e_2)=J(e_2)=-e_1$,
in such a way so as to adapt it to the underlying  Hermitian  metric
so that $\omega(x,y)=g(x,Jy)$. Therefore, define $J(e_3)=-e_4$
and $J(e_4)=e_3$.

\medskip
\begin{Remark}
The Lie algebra $\h_3(D)$ just described, is the Lie algebra
denoted by $D_{4,1/2}$ in \cite{Ovando}.
\end{Remark}

\medskip
{\bf Construction.}
Now, starting with the Frobenius-K\"ahler Lie algebra $\h_3(D)$,
we shall now produce a $5$-dimensional Sasakian Lie algebra
with trivial center. In order to do this, it suffices to choose a
derivation $E \in\Der(\h_3(D))$, satisfying
$E(e_1) = -e_2$, $E(e_2)=e_1$, and, $E(e_i)=0$ for $i=3,4$.
Observe that $\alpha \circ E = 0$, and that $[E,J]=0$. Therefore,
$\h_3(D,E):= (\h_3(D))(E)$ is a Sasakian Lie algebra
(see  {\bf Theorem \ref{FKtoSasaki}}).

\medskip
\begin{Remark}
The Lie algebra $\h_3(D,E)$ just described, is the Lie algebra
denoted by $\g_0$ in \cite{A-F-V}.
\end{Remark}

\medskip
On the other hand, since $\h_3(D)$ is a Frobenius-K\"ahler Lie algebra
with Frobenius structure given by $\alpha=e^3$, we may produce the
central extension of $\h_3(D)$ defined by the $2$-cocyle 
$\omega=-d\alpha$.
This construction yields again a Sasakian Lie algebra; this time, however, 
with non-trivial center.

\medskip
\begin{Remark}
The central extension of $\h_3(D)$ defined by the 2-cocycle
$\omega=-d\alpha$, is the Lie algebra denoted by 
$\g_5$ in  {\bf Theorem 10} of \cite{A-F-V}.
\end{Remark}

\medskip
Observe that if one starts with a Sasakian Lie algebra $\g$
defined by a given triple $(\xi,\alpha,\Phi)$, it seems easier to
construct first a Frobenius-K\"ahler Lie algebra $\g(D)$ through
a principal derivation $D\in\Der(\g)$ satisfying $[D,\Phi]=0$,
and then proceed to find a central extension of $\g(D)$
by means of the $2$-cocycle defined by its symplectic
form $\omega$. In this way one would be producing
a {\it reversed double extension\/,} in the sense that
one would adjoin the derivation $D$ first to obtain 
the Frobenius-K\"ahler Lie algebra $\g(D)$, and perform
the central extension by its $2$-cocyle afterwards.
That is, the operations are done exactly in the reverse
order as they were originally performed by V. Ka\c{c} to
produce his double extension procedure.

\begin{Remark}
It has been proved in \cite{B-R-S} that there are contact Lie algebras
that do not admit principal derivations. Thus, if one starts with such 
a Lie algebra that furthermore happens to be Sasakian,
then one cannot perform on it the {\it reversed double extension\/}
procedure to produce again a new (2 dimensions bigger) Sasakian Lie algebra.
\end{Remark}

\vskip1cm

\section*{Disclosure Statement}
The authors have no conflict of interest to declare that are relevant to this article.

\section*{Data Availability}
No datasets were generated or analyzed during the current study.

\section*{Acknowledgements} 
MCRV, GS and OASV would like to acknowledge the partial support received
by CONACyT grant A-S1-45886.
MCRV and GS would like to acknowledge the partial support received
by PRODEP grant UASLP-CA-228.
ASV also acknowledges the support from MB1411.

\bibliographystyle{plain}

\end{document}